\documentclass{amsart}
\usepackage{amsmath,amssymb,url}

\newcommand\F{{\mathbb F}}
\newcommand\Z{{\mathbb Z}}

\newcommand\SL{\mbox{\rm SL}}
\newcommand\GL{\mbox{\rm GL}}

\newcommand\gen[1]{\langle#1\rangle}
\newcommand\sz[1]{\left|#1\right|}

\author{Alexander Hulpke}
\address{Department of Mathematics, Colorado State University, 1874 Campus
Delivery, Fort Collins, CO, 80523-1874}
\email{hulpke@colostate.edu}

\thanks{The author's work has been supported in part by
NSF~Grant~DMS-1720146, and Simons Foundation Grant~852063,
which are gratefully acknowledged.}

\title{Proving infinite index for a subgroup of matrices}
\begin{document}
\begin{abstract}
We show how tools from computational group theory can be used to prove that
a subgroup of matrices has infinite index.
\end{abstract}
\maketitle
The combination of algorithmic methods for matrix
groups~\cite{baarnhielmholtleedhamobrien} with those for finitely presented
groups~\cite{simsbook} has been used
successfully~\cite{detinkoflanneryhulpkesp4} to prove that
certain subgroups of infinite matrix groups have finite index. The purpose of this note is to show, in
a concrete example, how the same toolkit can be used to prove the
infinite index of a particular subgroup. All calculations were performed in {\sf
GAP}~\cite{GAP4}, a transcript of the calculation and code 
is available at~\url{https://www.math.colostate.edu/~hulpke/paper/MaxIndexTranscript.txt}.

The group in question in this example is $\Gamma=\GL_2(\Z(\zeta))$, for $\zeta$ a
primitive 3rd root of unity. Motivated by~\cite{baechlemaheshwarymargolis}, 
A. B\"achle asked (private communication), whether the subgroup
$S=\gen{m_1,m_2,m_3,m_i,m_j,m_t}$, with
\begin{eqnarray*}
&&
%m_1=\left(\begin{array}{cc}
%97\zeta^2& -112\zeta-56\zeta^2\\
%112\zeta+56\zeta^2& 97\zeta^2
%\end{array}\right),
%\qquad
%m_i=\left(\begin{array}{cc}
%0& 1\\
%-1& 0
%\end{array}\right),\\
%&&
%m_2=\left(\begin{array}{cc}
%56\zeta+41\zeta^2& 56\zeta+112\zeta^2\\
%                56\zeta+112\zeta^2& -56\zeta+153\zeta^2
%\end{array}\right),
%\qquad
%m_j=\left(\begin{array}{cc}
%\zeta& \zeta^2\\
%\zeta^2& -\zeta
%\end{array}\right),\\
%&&
%m_3=\left(\begin{array}{cc}
%56\zeta+209\zeta^2& -56\zeta+56\zeta^2\\
%                -56\zeta+56\zeta^2& -56\zeta-15\zeta^2
%
m_1=\zeta
\left(\begin{array}{rr}
0& -112\\
112&0
\end{array}\right)+\zeta^2
\left(\begin{array}{rr}
97&-56\\
56&97
\end{array}\right)
,
\qquad
m_i=\left(\begin{array}{rr}
0& 1\\
-1& 0
\end{array}\right),\\
&&
m_2=56\zeta\left(\begin{array}{rr}
1&1\\
1&-1
\end{array}\right)+\zeta^2
\left(\begin{array}{rr}
41&112\\
112&153
\end{array}\right),
\qquad
m_j=\zeta\left(\begin{array}{rr}
1& \zeta\\
\zeta&-1
\end{array}\right),\\
&&
m_3=56\zeta\left(\begin{array}{rr}
1&-1\\
-1&-1
\end{array}\right)+\zeta^2
\left(\begin{array}{rr}
209&56\\
56&-15
\end{array}\right),
\qquad
m_t=\zeta\left(\begin{array}{cc}
1& 0\\
0&1
\end{array}\right)
\end{eqnarray*}
has finite index in $\Gamma$.

We obtain a finite presentation for $\Gamma$ from~\cite[Theorem~6.1]{swan71},
which gives us that $\Gamma\cong\gen{t,u,j,l,a,w}$ with
\begin{eqnarray*}
&&
t=\left(\begin{array}{cc}
1&1\\0&1
\end{array}\right),\qquad
u=\left(\begin{array}{cc}
1&\zeta\\0&1
\end{array}\right),\qquad
j=\left(\begin{array}{cc}
-1&0\\0&-1
\end{array}\right),\\
&&
l=\left(\begin{array}{cc}
\zeta^2&0\\0&\zeta
\end{array}\right),\qquad
a=\left(\begin{array}{cc}
0&-1\\1&0
\end{array}\right),\qquad
w=\left(\begin{array}{cc}
-\zeta&0\\0&1
\end{array}\right),
\end{eqnarray*}
subject to the relations/relators
\begin{eqnarray*}
&&
tu=ut, j^2, tj=jt, uj=ju, lj=jl , aj=ja ,l^3, l^{-1}tl=t^{-1}u^{-1},\\
&&
l^{-1}ul=t, a^2=j, (al)^2=j, (ta)^3=j, (ual)^3=j, wj=jw, w^6, wtw^{-1}=u^{-1},\\
&& wuw^{-1}=tu, waw^{-1}=jl^2a, wl=lw
\end{eqnarray*}
An application of Tietze transformations~\cite[\S 5.3.3]{holtbook}, 
eliminating the redundant generators
$u=wt^{-1}w^{-1}$, $j=a^2$, $l=w^{-1}a^{-1}wa^{-1}$, gives us
$\Gamma=\gen{t,a,w}$.

Using a norm-based reduction, as described in~\cite{hulpkespword}, we can obtain
expressions for the generators of $S$ as words in the generators
of $\Gamma$ as:
\begin{eqnarray*}
m_1&=&w*(tawt)^-8/w,\\
m_2&=&
w^{-1}aw^{-1}(t^{-1}w^{-1}t^{-1}a^{-1}twta^{-1})^3t^{-1}w^{-1}t^{-1}a^{-1}twt,\\
m_3&=& w^{-1}a^{-1}w^{-1}(a^{-1}twta^{-1}t^{-1}w^{-1}t^{-1})^4a^{-1},\\
m_i&=&a^{-1},\\
m_j&=&w^2a^{-1}t^{-1}w^{-1}a^{-1}t^{-1}/w\\
m_t&=&(w/a)^2.
\end{eqnarray*}
Attempts to determine the index if $S$ in $\Gamma$ by coset enumeration fail.

Following the approach of~\cite{detinkoflanneryhulpke17}, we next look at congruence images. Since
$\Gamma$ does not satisfy the congruence subgroup property this will only ever
provide a lower bound for the index of the subgroup,  even if it is finite.

Let $\varphi$ be the reduction modulo $2$ map
$\Gamma\to\GL_2(\F_4)$. We find that $\sz{\varphi(S)}=12$, while
$\sz{\varphi(\Gamma)}=180$, showing that $[\Gamma:S]\ge 15$ and in particular that
$\Gamma\not=S$.

Working simultaneously modulo $3,7,31,97,169,361,607$ (primes, found
experimentally, modulo which
a 3rd root of unity exists,
and modulo which the image of $S$  has
index at least $p/2$ in the respective $\GL_2(p)$) finds a quotient of
$\Gamma$
of order
\[
2^{43}3^{15}{5^3}{7^4}{13^5}{19^6}31{\cdot}97{\cdot}{101^2}607
\]
in which the image of $S$ has index
$2^{27}3^{10}{5^2}{7^3}13{\cdot}{19^2}31{\cdot}97{\cdot}{101^2}607$.

An index so large clearly puts coset enumeration outside the range of
feasibility. The fact that it has been so easy to build up so large an index
(indeed we could have tested further primes and would have obtained an even
larger index), however indicates that the subgroup might in fact be of
infinite index, and this is what we will show now:
\medskip

The method we shall use is to use a normal subgroup $N\lhd\Gamma$ of
finite index, that has an infinite abelian quotient, and so that $S\cap N$
has small index in $S$. (The former is possible only
because $\Gamma$ does not satisfy the congruence subgroup property.)
We then calculate generators for $N\cap S$ and show that $[N:N\cap S]$ is
infinite, contradicting that $S$ could have finite index in $\Gamma$.
\medskip

Concretely, let $\varphi$ be the reduction on $\Gamma$ modulo $7$. Let
$N=\ker\varphi$, then $[G:N]=\sz{\GL_2(7)}=2016$. We also calculate
$\sz{\varphi(S)}=24$.
This means that $N\cap S$ has index $\sz{\varphi(S)}=24$ in $S$.

We now construct generators for $N\cap S$.  $S$ acts, through $\varphi$ on
$\GL_2(7)$ by right multiplication, and the identity matrix has an orbit of
length $\sz{\varphi(S)}=24$. The stabilizer of the identity matrix is $N\cap
S$, and we find generators (121 of them) of $S\cap N$ as Schreier
generators.  

Next we use Reidemeister-Schreier rewriting~\cite[\S 5.3]{holtbook}
to obtain a presentation for $N$. Using this presentation, a Smith normal
form calculation~\cite[\S 9.2]{holtbook} will construct a homomorphism $\mu$ on
$N$ such that $\ker\mu=N'$.
In this example we find that $\mu(N)\cong\Z^8$. While \textsf{GAP} represents
the infinite factor group as a finitely presented group, a simplification
with Tietze transformations results in a presentation
on 8 generators. In this presentation the isomorphism to $\Z^8$ is given
simply by considering exponent sums of words.

We find that $\mu(N\cap S)$ is a submodule of $\Z^8$ of rank $3$. Thus
$[N:N\cap S]=\infty$ and we thus have shown that $[G:S]=\infty$, as claimed.
\medskip

Unfortunately, this approach will not work in cases when the congruence
subgroup property holds. We give a justification of this for the case of
$\Gamma=\SL_n(\Z)$, $n\ge 3$:
A subgroup of finite index in $\Gamma$ (our $N$ in the previous argument) is a
congruence subgroup, and thus the pre-image of a subgroup of a congruence
image $\varphi(\Gamma)$ for some congruence map $\varphi$, say the
congruence is modulo $m$. If $N$ has an infinite abelian quotient, we can
find quotients $\alpha_p$ of $N$ that are cyclic of any prime order $p$, in
particular for primes that are larger than any primes in the order of
$\varphi(\Gamma)$. But that means that for arbitrary large primes
$p$, $\Gamma$ has congruence images, whose order is a multiple of $p$,
without involving a composition factor $\mbox{PSL}_n(p)$, in contradiction
to the structure of congruence images of $\SL_n(\Z)$.

\bibliographystyle{amsplain}
\bibliography{mrabbrev,litprom}

\end{document}